# A non-balanced staggered-grid finite-difference scheme for the first-order elastic wave-equation modeling

Wenquan Liang[a]   Yanfei Wang[b,c,d], Ursula Iturrarán-Viveros[e]

[a]School of Resource Engineering, Longyan University, Longyan 364000, People's Republic of China

[b]Key Laboratory of Petroleum Resources Research, Institute of Geology and Geophysics, Chinese Academy of Sciences, Beijing 100029, People's Republic of China

[c]University of Chinese Academy of Sciences, Beijing 100049, People's Republic of China

[d]Institutions of Earth Science, Chinese Academy of Sciences, Beijing 100029, People's Republic of China

Corresponding author E-mail: yfwang@mail.iggcas.ac.cn

[e]Facultad de Ciencias, Universidad Nacional Autónoma de México, Circuito Escolar S/N, Coyoacán C.P. 04510; Ciudad de México

E-mail:ursula@ciencias.unam.mx

**ABSTRACT**

We introduce an efficient and accurate staggered-grid finite-difference (SGFD) method to solve the two-dimensional elastic wave equation. We use a coupled first-order stress-velocity formulation. In the standard implementation of SGFD method the same SGFD operator is used to approximate the spatial derivatives. However, we propose a numerical method based on mixed SGFD operators which happen to be more efficient with similar accuracy in comparison to uniform SGFD operator. We refer the proposed method as the non-balanced SGFD numerical scheme which means combining high-order SGFD operators with second-order SGFD operators. A very care attention is directed at the derivation of the SGFD operator coefficients. The correctness of proposed scheme is proven by dispersion analysis. Through SGFD modeling examples, we verify/demonstrate that the proposed non-balanced operator offers a similar level of accuracy with a cheaper computation cost compared to the more expensive balanced SGFD method.

**INTRODUCTION**

The wave propagation in elastic solid has practical applications in seismic exploration, seismology, material science, geotechnical engineering and other areas of research. To numerically simulate elastic wave propagation, various numerical methods have been developed. Examples are the finite-difference, finite-volume or finite-element method. The probably most common method is the finite-difference method. In the finite- difference method the spatial derivatives of partial differential equation is approximated by finite-difference on an equally spaced grid. The strength of finite-difference is its simple numerical representation. The finite-difference method is often accompanied by numerical dispersion and predominantly observed in the surface wave simulation.

To suppress the issue of numerical dispersion an approach of staggered grid finite-difference method (SGFD) is effectively employed and rigorously studied by various researchers (Virieux 1984; Levander, 1988; Saenger et al, 2000; Kantartzis et al, 2006; Moczo et al, 2002, 2007; Lisitsa and Vishnevsky 2010; Sun et al, 2012; Bohlen et al. 2016; Etemadsaeed et al. 2016; Igel 2017). In SGFD method material properties are distributed in staggered fashion at grid points which smoothened out a property contrast across the interfaces. The SGFD method inherits the properties a generic finite-difference method naturally has i.e. simplicity and efficiency.

The accuracy of the SGFD is primarily controlled by the discretization parameters i.e. time and space step. However, decreasing the spacing (time or space) to achieve the high order accuracy affects the efficiency of the scheme. As an alternative, coarse-grid methods were developed in the form of high-order finite-difference methods, where large space steps are used in combination with long finite-difference stencils that provides with a reasonable accuracy (Kelly, 1976; Holberg, 1987). This high-order finite-difference method is now the prevalently adopted numerical scheme to obtain the high accuracy. For example, a finite-difference operator at grid point using 21 grid points (also called as a stencil size) will result in an accurate approximation of second-order spatial derivative (Stork, 2013). The large stencils typically provide with a higher accuracy and reduce the number of grid points per minimum wavelength. However, implementation of the long stencils results into an increased computational cost and especially associated with the memory requirements.

In this paper, we propose a numerical method that reduces the cost of the coarse grid methods using the large stencils, by using two different finite-difference operators (a long and a short one) of different size and used to compute the spatial derivatives. Both the long and the short stencils are used simultaneously to evaluate the spatial derivatives within a time step. A careful approach is considered to compute the SGFD coefficients for long stencil, which cancels out the error produced by short operators. We refer to this as a 'non-balanced' SGFD scheme.

The approach of using the finite-difference operators of different sizes in a single simulation has been already proposed by other researchers (e.g., Stork 2013; Kindelan et al. 2016), but in this work, we combines a long finite- difference stencil with only a short second-order finite-difference operator which is also a novelty of the current work. A similar method has already been successfully applied to solve the acoustic wave equation (Liang et al., 2018a, b). In this paper we illustrate the efficiency of this method by applying it to the elastic wave equation. The proposed method achieves accuracy similar to that of coarse-grid and long stencil methods but with 40% less computational cost.

The paper is structured as follows: First, we rigorously define the 'non-balanced' SGFD method. Then, we define the method to compute the finite-difference coefficients using the dispersion relation of the implementation. Then, we analyze the stability of the proposed numerical scheme. Finally, we perform the computational experiments for two different elastic models of varying complexity and show the efficacy of the proposed 'non-balanced' SGFD scheme.

## THE NON-BALANCED SGFD SCHEME

The time domain equations for wave propagation in a heterogeneous elastic medium can be found in (Virieux, 1984). The two-dimensional (2D) velocity-stress equation for elastic wave propagation in the (x, z) plane is expressed as

i) Euler-Newton equations:

$$\frac{\partial v_x}{\partial t} = \frac{1}{\rho}\left(\frac{\partial \tau_{xx}}{\partial x} + \frac{\partial \tau_{xz}}{\partial z}\right), \tag{1}$$

$$\frac{\partial v_z}{\partial t} = \frac{1}{\rho}\left(\frac{\partial \tau_{xz}}{\partial x} + \frac{\partial \tau_{zz}}{\partial z}\right), \tag{2}$$

ii) Constitutive equations:

$$\frac{\partial \tau_{xx}}{\partial t} = (\lambda + 2\mu)\frac{\partial v_x}{\partial x} + \lambda \frac{\partial v_z}{\partial z}, \tag{3}$$

$$\frac{\partial \tau_{zz}}{\partial t} = \lambda \frac{\partial v_x}{\partial x} + (\lambda + 2\mu)\frac{\partial v_z}{\partial z}, \tag{4}$$

$$\frac{\partial \tau_{xz}}{\partial t} = \mu(\frac{\partial v_x}{\partial z} + \frac{\partial v_z}{\partial x}), \tag{5}$$

where $(\tau_{xx}, \tau_{zz}, \tau_{xz})$ represents the stress components, $(v_x, v_z)$ is the velocity vector, $\lambda$ and $\mu$ are Lame coefficients and $\rho$ is the density of the medium.

To solve the system of elastic wave equation (1)-(5) numerically, Virieux (Virieux, 1984) discretized the spatial partial derivatives in equations (1)-(5) using the SGFD method. Virieux (Virieux, 1984) used the same SGFD operator to approximate the spatial derivatives appearing in equations (1)-(5). We refer the SGFD scheme used by Virieux (Virieux, 1984) as the *balanced SGFD scheme*.

We introduce different SGFD operators, refereed as non-balanced SGFD scheme, for different spatial derivatives in equations (1)-(5). The SGFD numerical schemes for equations (1) –(5) are expressed as

$$\frac{v_x\begin{bmatrix}1\\0,0\end{bmatrix}-v_x\begin{bmatrix}0\\0,0\end{bmatrix}}{\Delta t}=\frac{1}{\rho}\left\{\frac{1}{h}\sum_{m=1}^{M}c_m\left[\tau_{xx}\begin{bmatrix}1/2\\m-1/2,0\end{bmatrix}-\tau_{xx}\begin{bmatrix}1/2\\-m+1/2,0\end{bmatrix}\right]+\frac{1}{h}\left[\tau_{xz}\begin{bmatrix}1/2\\0,1/2\end{bmatrix}-\tau_{xz}\begin{bmatrix}1/2\\0,-1/2\end{bmatrix}\right]\right\}, (6)$$

$$\frac{v_z\begin{bmatrix}1\\1/2,1/2\end{bmatrix}-v_z\begin{bmatrix}0\\1/2,1/2\end{bmatrix}}{\Delta t}=\frac{1}{\rho}\left\{\frac{1}{h}\left[\tau_{xz}\begin{bmatrix}1/2\\1,1/2\end{bmatrix}-\tau_{xz}\begin{bmatrix}1/2\\0,1/2\end{bmatrix}\right]+\frac{1}{h}\sum_{m=1}^{M}c_m\left[\tau_{zz}\begin{bmatrix}1/2\\1/2,m\end{bmatrix}-\tau_{zz}\begin{bmatrix}1/2\\1/2,-m+1\end{bmatrix}\right]\right\}, (7)$$

$$\frac{\tau_{xx}\begin{bmatrix}1/2\\1/2,0\end{bmatrix}-\tau_{xx}\begin{bmatrix}-1/2\\1/2,0\end{bmatrix}}{\Delta t}=\frac{(\lambda+2\mu)}{h}\left[v_x\begin{bmatrix}0\\1,0\end{bmatrix}-v_x\begin{bmatrix}0\\0,0\end{bmatrix}\right]+\frac{\lambda}{h}\left[v_z\begin{bmatrix}0\\1/2,1/2\end{bmatrix}-v_z\begin{bmatrix}0\\1/2,-1/2\end{bmatrix}\right], (8)$$

$$\frac{\tau_{zz}\begin{bmatrix}1/2\\1/2,0\end{bmatrix}-\tau_{zz}\begin{bmatrix}-1/2\\1/2,0\end{bmatrix}}{\Delta t}=\frac{\lambda}{h}\left[v_x\begin{bmatrix}0\\1,0\end{bmatrix}-v_x\begin{bmatrix}0\\0,0\end{bmatrix}\right]+\frac{(\lambda+2\mu)}{h}\left[v_z\begin{bmatrix}0\\1/2,1/2\end{bmatrix}-v_z\begin{bmatrix}0\\1/2,-1/2\end{bmatrix}\right], (9)$$

$$\frac{\tau_{xz}\begin{bmatrix}1/2\\0,1/2\end{bmatrix}-\tau_{xz}\begin{bmatrix}-1/2\\0,1/2\end{bmatrix}}{\Delta t}=\frac{\mu}{h}\sum_{m=1}^{M}c_m\left[v_z\begin{bmatrix}0\\m-1/2,1/2\end{bmatrix}-v_z\begin{bmatrix}0\\-m+1/2,1/2\end{bmatrix}+v_x\begin{bmatrix}0\\0,m\end{bmatrix}-v_x\begin{bmatrix}0\\0,-m+1\end{bmatrix}\right], (10)$$

where $u\begin{bmatrix}k\\i,j\end{bmatrix}=u(x+ih, z+jh, t+k\Delta t)$, $u=[v_x, v_z, \tau_{xx}, \tau_{zz}, \tau_{xz},]$, $M$ is the SGFD operator length, $c_m$ (for each *m)* is the SGFD operator coefficient, $\Delta t$ is the temporal grid interval and *h* is the distance between the two consecutive space grid points.

The key component in proposed non-balanced SGFD scheme is the usage of the short SGFD operator for half of the spatial derivatives in equations (6)-(10). The numerical schemes in equations (6)-(10) are defined as the *non-balanced SGFD scheme* as it is not balanced for its SGFD operator length.

At first, the non-balanced SGFD scheme equations (6)-(10) may not seem attractive as spatial derivatives approximated with second-order SGFD operators would results into spatial grid dispersions. However, specially designed SGFD coefficients for the non-balanced SGFD scheme can result into reduced grid dispersion and improved accuracy of the numerical simulation.

**COMPUTATION OF COEFFICIENTS OF THE NON-BALANCED SGFD OPERATOR**

Substitution of equations (3)-(5) into equations (1)-(2) yields second-order elastic wave equation. The matrix form of the second-order elastic wave-equation in the frequency-wavenumber domain is expressed as

$$\begin{pmatrix} -\alpha^2 k_x^2 - \beta^2 k_z^2 + \omega^2 & -(\alpha^2 - \beta^2) k_x k_z \\ -(\alpha^2 - \beta^2) k_x k_z & -\alpha^2 k_z^2 - \beta^2 k_x^2 + \omega^2 \end{pmatrix} \begin{pmatrix} v_x \\ v_z \end{pmatrix} = 0, \qquad (11)$$

where $k_x$ and $k_z$ are the spatial wave numbers in the x and z directions respectively, $\beta$ and $\alpha$ are the S- and P- wave velocities respectively with

$$\alpha^2 = (\lambda + 2\mu)/\rho, \tag{12}$$

$$\beta^2 = \mu/\rho. \tag{13}$$

We show why the non-balanced SGFD scheme and the balanced SGFD scheme are equivalent for the elastic wave-equation. Equation (11) can be written as

$$\left[-\alpha^2 k_x^2 - \beta^2 k_z^2 + \omega^2\right]\left[-\alpha^2 k_z^2 - \beta^2 k_x^2 + \omega^2\right] - \left[\left(\alpha^2 - \beta^2\right)k_x k_z\right]^2 = 0. \tag{14}$$

For the balanced SGFD scheme

$$k_x = \frac{2}{h}\sum_{m=1}^{M} c_m \left(\sin\left[(m-0.5)k_x h\right]\right), \tag{15}$$

$$k_z = \frac{2}{h}\sum_{m=1}^{M} c_m \left(\sin\left[(m-0.5)k_z h\right]\right), \tag{16}$$

$$k_x k_z = \left[\frac{2}{h}\sum_{m=1}^{M} c_m \left(\sin\left[(m-0.5)k_x h\right]\right)\right]\left[\frac{2}{h}\sum_{m=1}^{M} c_m \left(\sin\left[(m-0.5)k_z h\right]\right)\right]. \tag{17}$$

For the non-balanced SGFD scheme

$$k_x^2 = \frac{2}{h}\sin\left[(0.5)k_x h\right]\left[\frac{2}{h}\sum_{m=1}^{M} c_m \left(\sin\left[(m-0.5)k_x h\right]\right)\right], \tag{18}$$

$$k_z^2 = \frac{2}{h}\sin\left[(0.5)k_z h\right]\left[\frac{2}{h}\sum_{m=1}^{M} c_m \left(\sin\left[(m-0.5)k_z h\right]\right)\right], \tag{19}$$

and $k_x k_z$ has two different SGFD operator approximations

$$k_x k_z = \frac{4\sin(0.5 k_x h)}{h^2}\left[\sum_{m=1}^{M} c_m \left(\sin\left[(m-0.5)k_z h\right]\right)\right], \tag{20}$$

or

$$k_x k_z = \frac{4\sin(0.5 k_z h)}{h^2}\left[\sum_{m=1}^{M} c_m \left(\sin\left[(m-0.5)k_x h\right]\right)\right]. \tag{21}$$

We have demonstrated that the non-balanced SGFD scheme can approximate $k_x^2, k_z^2$ good as the balanced SGFD scheme in our previous work. This means (Liang et al, 2018a, b)

$$k_x^2 h^2 \approx \left[\sum_{m=1}^{M} 2c_m \left(\sin\left[(m-0.5)k_x h\right]\right)\right]^2 \approx 2\sin\left[(0.5)k_x h\right]\left[2\sum_{m=1}^{M} c_m \left(\sin\left[(m-0.5)k_x h\right]\right)\right], \quad (22)$$

$$k_z^2 h^2 \approx \left[\sum_{m=1}^{M} 2c_m \left(\sin\left[(m-0.5)k_z h\right]\right)\right]^2 \approx 2\sin\left[(0.5)k_z h\right]\left[2\sum_{m=1}^{M} c_m \left(\sin\left[(m-0.5)k_z h\right]\right)\right]. \quad (23)$$

Neither equation (20) nor (21) is equivalent to equation (17). However, it is $k_x^2 k_z^2$ in equation (14). From equation (17) and equation (20), (21) we can get

$$\left(k_x k_z h^2\right)^2 \approx \left[2\sum_{m=1}^{M} c_m \left(\sin\left[(m-0.5)k_x h\right]\right)\right]^2 \left[2\sum_{m=1}^{M} c_m \left(\sin\left[(m-0.5)k_z h\right]\right)\right]^2 \approx$$
$$2\sin(0.5k_x h) 2\sin(0.5k_z h)\left[2\sum_{m=1}^{M} c_m \left(\sin\left[(m-0.5)k_x h\right]\right)\right]\left[2\sum_{m=1}^{M} c_m \left(\sin\left[(m-0.5)k_z h\right]\right)\right], \quad (24)$$

Equation (24) means the non-balanced SGFD scheme can approximate the mixed spatial derivative as good as the balanced SGFD scheme. Equation (22), (23) and (24) means the non-balanced SGFD scheme can approximate equation (14) as good as the balanced SGFD scheme. In a word, the balanced SGFD scheme tries to approximate $k_x, k_z$ and $k_x k_z$, while the proposed non-balanced SGFD scheme tries to approximate $k_x^2, k_z^2$ and $[k_x k_z]^2$.

In the following, we demonstrate how to get the SGFD coefficients for the non-balanced SGFD scheme. The two roots of equation (11) yield the following dispersion relation

$$\omega^2 = \frac{1}{2}\left(\alpha^2 + \beta^2\right)\left(k_x^2 + k_z^2\right) \pm \frac{1}{2}\left(\alpha^2 - \beta^2\right)\left(k_x^2 + k_z^2\right). \quad (25)$$

Then we would have

$$\alpha^2 k_x^2 + \alpha^2 k_z^2 - 4\sin^2\left(\omega\Delta t/2\right) = 0, \qquad (26)$$

$$\beta^2 k_x^2 + \beta^2 k_z^2 - 4\sin^2\left(\omega\Delta t/2\right) = 0. \qquad (27)$$

The SGFD operator coefficients could be designed either using equation (26) or (27). It is well known that the smaller wave propagation speed has severe grid dispersion. Consequently, equation (27) is used to determine the SGFD operator coefficient as $\beta < \alpha$. After substituting equation (15) and (16) into equation (27), the dispersion relation for the balanced SGFD is expressed as (Wang et al, 2014)

$$\left[\sum_{m=1}^{M} c_m \sin((m-0.5)k_x h)\right]^2 + \left[\sum_{m=1}^{M} c_m \sin((m-0.5)k_z h)\right]^2 = \frac{1-\cos(k\beta\Delta t)}{2r^2}, \qquad (28)$$

where $r = \beta\Delta t/h$ and $(k_x, k_z) = k(\cos\theta, \sin\theta)$, in which $\theta$ is the plane wave-propagation angle measured from the horizontal plane perpendicular to the z-axis.

Similarly, the dispersion equation for the non-balanced SGFD scheme is expressed as

$$\sin(0.5k_z h)\sum_{m=1}^{M} c_m \left(\sin(m-0.5)k_z h\right) + \sin(0.5k_x h)\sum_{m=1}^{M} c_m \left(\sin(m-0.5)k_x h\right) = \frac{1-\cos(k\beta\Delta t)}{2r^2}, \qquad (29)$$

where $c_m$ (for all $m$) are the SGFD operator coefficients.

It should be noted that, when determining the SGFD coefficients in the time-space domain using equation (29), the SGFD coefficients correspond to the wave propagation velocity $\beta$. Therefore, there are many SGFD coefficients satisfying the time-space domain methods.

Considering only one direction and taking the limit on $\Delta t$ as the time step $\Delta t$ approaches to zero, the equation (29) is expressed as

$$sin(0.5kh)\sum_{m=1}^{M} c_m \left( sin(m-0.5)kh \right) = \lim_{\Delta t \to 0} \frac{1-cos(k\beta\Delta t)}{2(\beta\Delta t/h)^2} = \frac{k^2h^2}{4}, \qquad (30)$$

which is exactly the same as in our previous work (Liang et al, 2018b). Therefore, with equations (29) and (30), we can calculate the SGFD coefficients for the non-balanced SGFD scheme in the both time and space domain or in the space domain only.

Readers are advised to refer to our previous works (Liang et al, 2018a, 2018b) for detailed explanation on computation of SGFD coefficients for the non-balanced SGFD scheme using equations (29) and (30). By using this approach, we obtain the non-balanced SGFD coefficients with a minimal L2-norm error for equations (29)-(30). These coefficients are valid over a large wavenumber range for a number of directions.

Thus, one can get the optimized SGFD coefficient in the space domain (with equation 15 or 16) or in the time-space domain (with equation 28) for the balanced SGFD scheme. However, optimized non-balanced SGFD coefficient could be obtained in the time-space and space domain by using equation (29) and (30), respectively.

As equation (30) is independent of the velocity, thus there will be only one group of SGFD coefficients. Additionally, the time dispersion can be eliminated by using the method proposed by Koene et al. (2017).

The commonly used non-balanced first-order SGFD coefficients computed in the space domain for equations (6)-(10) are shown in Table 1 (obtained by the least-squares solution of equation 30). These coefficients are also used for the computational experiment presented in this paper.

DISPERSION ANALYSES

We show a comparison between the numerical errors in solutions obtained from the balanced and non-balanced SGFD scheme, in the time-space domain. The numerical error $\delta$ of the balanced SGFD scheme for the $P$ wave is described by (Liu and Sen, 2011)

$$\delta(\theta) = \frac{\alpha_{FD}}{\alpha} = \frac{2}{rkh} \sin^{-1}\left(r\sqrt{q_1}\right), \tag{31}$$

where

$$q_1 = \left(\sum_{m=1}^{M} c_m \sin\left[(m-0.5)kh\cos\theta\right]\right)^2 + \left(\sum_{m=1}^{M} c_m \sin\left[(m-0.5)kh\sin\theta\right]\right)^2. \tag{32}$$

Similarly, the numerical error $\delta$ of the balanced SGFD scheme for the $S$ wave could be also obtained.

The numerical error $\delta$ of the non-balanced SGFD scheme described in this paper for the $P$ wave is defined by

$$\delta(\theta) = \frac{\alpha_{FD}}{\alpha} = \frac{2}{rkh}\sin^{-1}\left(r\sqrt{q_2}\right), \tag{33}$$

where

$$\begin{aligned} q_2 = &\sin(0.5kh\cos\theta)\sum_{m=1}^{M}c_m\sin\left[(m-0.5)kh\cos\theta\right] \\ &+ \sin(0.5kh\sin\theta)\sum_{m=1}^{M}c_m\sin\left[(m-0.5)kh\sin\theta\right]. \end{aligned} \tag{34}$$

Following the equation (33), the numerical error $\delta$ of the non-balanced SGFD scheme (described in this paper) for the $S$ wave could be also obtained.

Figures 1-3 illustrate the plots of numerical error of the different SGFD schemes for a homogeneous 2D elastic model. In all the figures, $h$ is the spatial grid interval, $M$ is the length of the SGFD operator and $\Delta t$ is the time step. Results shown in the Figure 1 is obtained by using the balanced SGFD scheme with the coefficients calculated by the Taylor expansion. The absolute value of the numerical error is less than 0.005 (within 70 % of the $kh$ range). Figure 2 represents error plots for the balanced SGFD scheme with the coefficients calculated by the least-squares method (we refer to it as the *optimized solution*). When we consider an error under 0.005, we note that the optimized SGFD operator coefficients can suppress the numerical error in a larger wavenumber range (Figure 2) than the non-optimized one (Figure 1), obtained from the Taylor expansion. This shows an improvement of the newly proposed non-balanced optimized SGFD scheme over non-optimized SGFD scheme. Figure 3 is computed using the proposed non-balanced SGFD scheme with the coefficients calculated by the least-squares method. From Figures 2 and 3, we draw the conclusion that the precision of the proposed non-balanced SGFD scheme is similar to the balanced one, whereas our proposed method needs less computational time due to the short operator length of equations (6)-(10).

We also show the dispersion error of the mixed derivative. For the balanced SGFD scheme,

$$E = \left(k_x k_z h^2\right)^2 - \left[2\sum_{m=1}^{M} c_m \left(\sin\left[(m-0.5)k_x h\right]\right)\right]^2 \left[2\sum_{m=1}^{M} c_m \left(\sin\left[(m-0.5)k_z h\right]\right)\right]^2. \tag{35}$$

For the non-balanced SGFD scheme,

$$E = \left(k_x k_z h^2\right)^2 - 16\sin(0.5k_x h)\sin(0.5k_z h)\sum_{m=1}^{M} c_m \left(\sin\left[(m-0.5)k_x h\right]\right)\sum_{m=1}^{M} c_m \left(\sin\left[(m-0.5)k_z h\right]\right). \tag{36}$$

From Figure 4, we can observe that the dispersion error of mixed derivative from the balanced SGFD scheme and the non-balanced SGFD scheme are almost the same.

## STABILITY ANALYSIS

From equations (26) and (27), we can obtain the following equations for the balanced SGFD scheme (Levander, 1988)

$$\sin^2(0.5\omega\tau) = \left(\frac{\Delta t \beta}{h}\right)^2 \left(X^2 + Z^2\right), \tag{37}$$

$$\sin^2(0.5\omega\tau) = \left(\frac{\Delta t \alpha}{h}\right)^2 \left(X^2 + Z^2\right), \tag{38}$$

where

$$X = \sum_{m=1}^{M} c_m \sin\left((m-0.5)k_x h\right), \tag{39}$$

$$Z = \sum_{m=1}^{M} c_m \sin\left((m-0.5)k_z h\right). \tag{40}$$

If both types of waves (*P* and *S*) are generated and propagated in a medium, the condition in equation (38) for the *P* wave has to be considered as the joint stability condition, due to $\alpha$ being large. Once the stability condition for $\alpha$ is satisfied, then the condition naturally extends for $\beta$.

From equation (38) we obtain

$$\left(\frac{\Delta t \alpha}{h}\right)^2 \left(X^2 + Z^2\right) \leq 1. \tag{41}$$

Equation 41 gives a stability condition for the balanced SGFD scheme (with $k_x h=\pi$, $k_z h=\pi$)

$$\left(\frac{\Delta t \alpha}{h}\right) \leq 1/\sqrt{\left(X^2 + Z^2\right)} \leq \frac{1}{\sqrt{2}\sum_{m=1}^{M}|c_m|}. \tag{42}$$

For the proposed non-balanced SGFD scheme, we get

$$X'^2 = \sin(0.5 k_x h) \sum_{m=1}^{M} c_m \left(\sin\left((m-0.5)k_x h\right)\right), \tag{43}$$

$$Z'^2 = \sin(0.5 k_z h) \sum_{m=1}^{M} c_m \left(\sin\left((m-0.5)k_z h\right)\right). \tag{44}$$

Similarly, the stability condition for the proposed non-balanced SGFD scheme is expressed as

$$r = \left(\frac{\Delta t \alpha}{h}\right) \leq 1/\sqrt{\left(X'^2 + Z'^2\right)} \leq s = \frac{1}{\sqrt{2\sum_{m=1}^{M}|c_m|}}. \tag{45}$$

We demonstrate the effectiveness of the stability condition in equation (45) with an example. We use the SGFD coefficients in the last row of Table 1. Let $h$=10m, $\Delta t$ =0.001s. The stability condition is calculated as $r \leq 1/\left(\sqrt{2\sum_{m=1}^{M}|c_m|}\right) = 0.491$. Therefore, from equation (45), if $\alpha$=4912 m/s, it should be unstable. However, if $\alpha$=4908 m/s, the scheme will be stable. We use the numerical simulation to verity the above stability condition with a homogeneous model. Our calculations show that if we set $\alpha = 4912$ m/s, $h = 10$ m and $\Delta t = 0.001$ s, the numerical process is unstable, whereas for $\alpha = 4908$ m/s, $h = 10$ m and $\Delta t = 0.001$ s, it will be stable.

## NUMERICAL SIMULATIONS

**Homogeneous model**

To validate the correctness of numerical scheme, we compare the numerical solution obtained for a homogeneous model. The discretization parameters for the SGFD operator are $h$=10m, $\Delta t$ =1 ms and $M$ = 7, $\alpha$ =1732.1 m/s and $\beta$ =1000 m/s. We consider a seismic source and receiver locations as shown in Figure 5. A Ricker wavelet with a central frequency at 14Hz is used as the seismic source.

The time history of particle velocities $v_z$ and $v_x$ are computed at 5 different receiver locations using the balanced SGFD scheme and the non-balanced SGFD scheme are shown in Figure 6. Figure 6(a)-(b) shows that seismograms are almost identical for the balanced and the non-balanced SGFD scheme with $M$=4[1]. There is space grid dispersion on both of them (as shown in the red eclipse of Figure 6a and b). Similarly, Figure 6(c)-(d) shows that seismograms are almost identical for the balanced and the non-balanced SGFD scheme with *M*=7. The space grid dispersion is suppressed with the increase of the SGFD operator length (as shown in the red eclipse of Figure 6c and d).Figure 6e, 6f, 6g,6h show the seismograms in figure 6c and 6d. From figure 6e-h, we can observe that the simulation results are almost identical for the balanced and the non-balanced SGFD scheme.

---

[1] All the figures in this paper are reproducible with code shared at Github
https://github.com/LongyanU/Geo-elastic3

**Salt velocity model**

Now, we consider the famous salt velocity model shown in Figure 7. The numerical simulation parameters for the SGFD operator are given as: *h*=10m, $\Delta t$ =1ms and M=7. In this example, we use long SGFD operator as the analytical method (M=30).

Figure 8(a) and 8(b) represents the seismic record of the $v_z$ component computed by using the balanced SGFD scheme (with non-optimized SGFD coefficient) and the non-balanced SGFD scheme (with optimized SGFD coefficient in the last row of Table 1), respectively. Figure 8(c) shows the seismic record of the $v_z$ component obtained by the balanced SGFD scheme with *M*=30 (with non-optimized SGFD coefficient). In Figure 8(d), the slices of seismograms are shown which are extracted from Figure 8(a)-(c). From Figure 8(d) and 8 (e), we observe that the results of the non-balanced SGFD scheme are very close to the balanced SGFD scheme with a very long SGFD operator (*M*=30).

Figure 9(a) and Figure 9(b) represents the seismic record of the $\tau_{xx}$ component computed with the balanced SGFD and non-balanced SGFD scheme, respectively.

Figure 9(c) shows the seismic record of the $\tau_{xx}$ component obtained by the balanced SGFD scheme with $M$=30. Figure 9(d) shows slices of seismograms extracted from Figure 9(a)-(c). From Figure 9(d) and 9(e), we may say that the results obtained from the non-balanced SGFD scheme and the balanced SGFD scheme with a very long SGFD operator ($M$=30) are in very good agreement.

In theory, assuming that the simulation time for the traditional balanced SGFD scheme is $M+M$, then the simulation time for the non-balanced SGFD scheme should be $M+1$. It takes 962 seconds to get the simulated wave field shown in Figure 8(a), whereas it takes only 553 seconds to get the simulated wave field shown in Figure 8(b). This concords with the analysis as 553/962 is equal to $(M+1)/(M+M)$. Therefore, the non-balanced SGFD scheme can save about 42% of simulation time compared with the traditional balanced SGFD scheme.

## CONCLUSION

In this paper, we first propose a non-balanced SGFD scheme for solving the 2D elastic wave equation. Then we present a rigorous method to determine the SGFD operator coefficients for the proposed non-balanced SGFD scheme. The numerical error in solution obtained from the non-balanced SGFD scheme is small when comparing it against balanced SGFD operator. Computational experiments are performed to illustrate the advantages of the proposed SGFD scheme. The proposed numerical scheme is also carried out for a salt velocity model and results are compared against existing implementation of SGFD schemes. The results show that with the non-balanced SGFD scheme and the corresponding FD operator coefficients, we require less computational time to compute the numerical simulation without losing accuracy. Therefore, we concluded that the non-balanced SGFD can be considered as a very useful tool to simulate the elastic wave-equations efficiently and accurately.

**LIST OF FIGURE CAPTIONS**

Figure 1. Numerical error curves of the balanced SGFD scheme with the non-optimized operator coefficients for 2D elastic wave equation with $\alpha$=2598m/s, $\beta$=1500m/s, $h$=20m, $M$=7, $\Delta t$ =1 ms: (a) $S$ wave and (b) $P$ wave.

Figure 2. Numerical error curves of the balanced SGFD scheme with the optimized operator coefficient for 2D elastic wave equation with $\alpha$=2598m/s, $\beta$=1500m/s, $h$=20m, $M$=7, $\Delta t$ =1ms: (a) $S$ wave and (b) $P$ wave.

Figure 3. Numerical error curves of the non-balanced SGFD scheme with the optimized operator coefficient for 2D elastic wave equation with $\alpha$=2598m/s, $\beta$=1500m/s, $h$=20 m, $M$=7, $\Delta t$ =1ms: (a) $S$ wave and (b) $P$ wave.

Figure 4. Dispersion error of the mixed derivative $\left(k_x k_z h^2\right)^2$ (a) The balanced SGFD scheme with $M$=7 and (b) the non-balanced SGFD scheme with $M$=7.

Figure 5. Source and receiver location

Figure 6. Comparison of seismograms obtained at the 1,3,5,7,9 receivers from the top to the bottom (a) the $v_z$ component obtained with the balanced and the non-balanced SGFD scheme with $M$=4; (b) the $v_x$ component obtained with the balanced and the non-balanced SGFD scheme with $M$=4; (c) the $v_z$ component obtained with the balanced and the non-balanced SGFD scheme with $M$=7; (d) the $v_x$ component obtained with the balanced and the non-balanced SGFD scheme with $M$=7;(e) the local enlargement of the first seismogram in figure (c); (e) the local enlargement of the first seismogram in figure (d); (g) the local enlargement of the third seismogram in figure (c); (h) the local enlargement of the third seismogram in figure (d);

Figure 7. Salt velocity model: (a) *P* wave velocity model; (b) *S* wave velocity model.

Figure 8. Seismic records of the $v_z$ component: (a) the balanced SGFD scheme (M=7); (b) the non-balanced SGFD scheme (M=7); (c) the balanced SGFD scheme (M=30); (d) slices of the seismograms extracted from (a)-(c) at x/dx=222; (e) difference of seismograms from (d).

Figure 9. Seismic records of $\tau_{xx}$ component: (a) the balanced SGFD scheme; (b) the non-balanced SGFD scheme; (c) the balanced SGFD scheme (M=30); (d) slices of seismograms extracted from (a)-(c) at x/dx=222; (e) difference of seismograms from (d).

## LIST OF TABLES

**Table 1.** Optimized first-order SGFD coefficients for the non-balanced FD scheme.

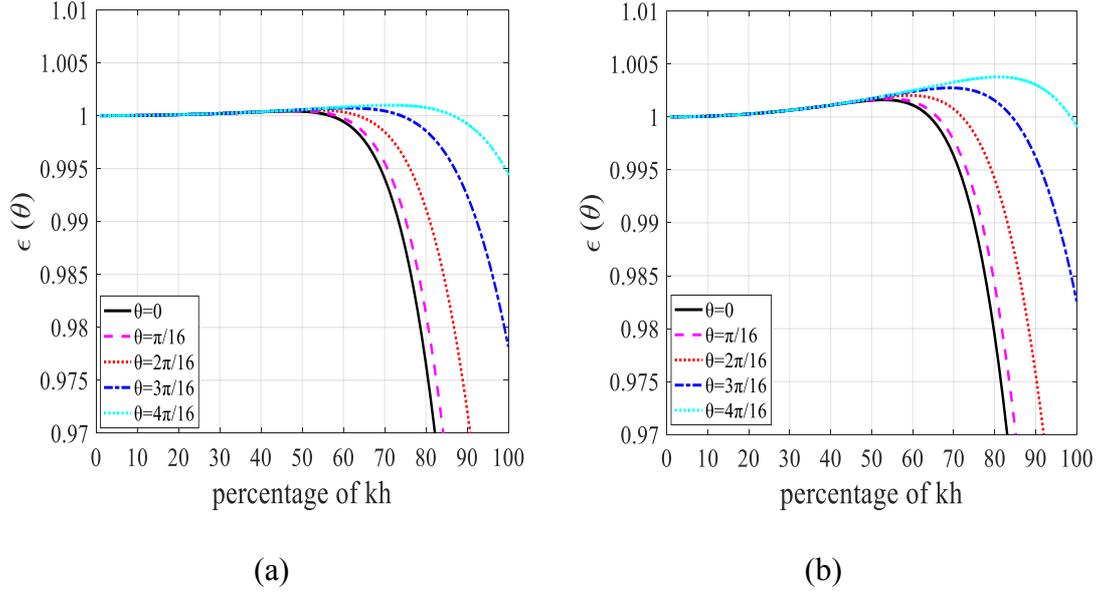

(a)                  (b)

Figure 1. Numerical error curves of the balanced SGFD scheme with the non-optimized operator coefficients for the 2D elastic wave equation with $\alpha$=2598m/s, $\beta$=1500m/s, $h$=20m, $M$=7, $\Delta t$ =1 ms. (a) $S$ wave and (b) $P$ wave.

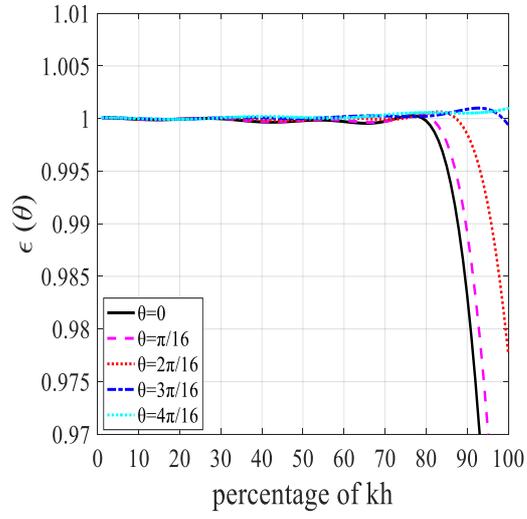 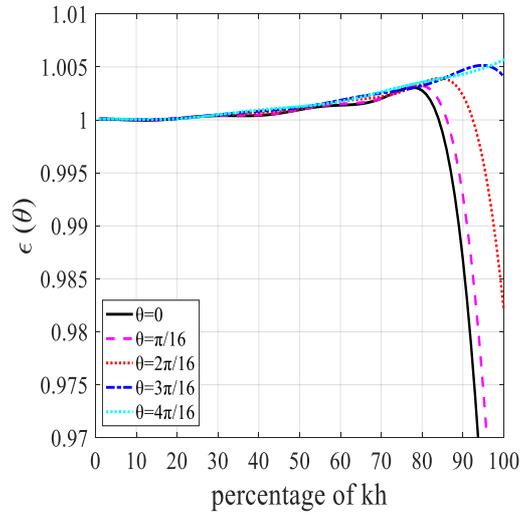

(a)                  (b)

Figure 2. Numerical error curves of the balanced SGFD scheme with the optimized operator coefficient for the 2D elastic wave equation with $\alpha$=2598m/s, $\beta$=1500m/s, $h$=20m, $M$=7, $\Delta t$ =1ms. (a) $S$ wave and (b) $P$ wave.

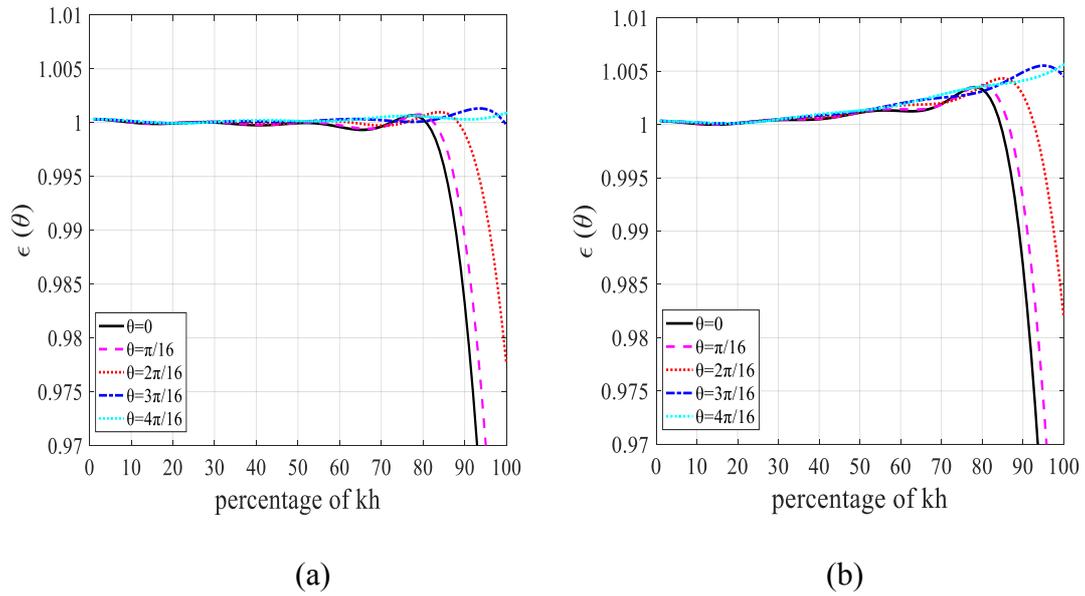

Figure 3. Numerical error curves of the non-balanced SGFD scheme with the optimized operator coefficient for the 2D elastic wave equation with $\alpha$=2598m/s, $\beta$=1500m/s, $h$=20 m, $M$=7, $\Delta t$ =1ms. (a) $S$ wave and (b) $P$ wave.

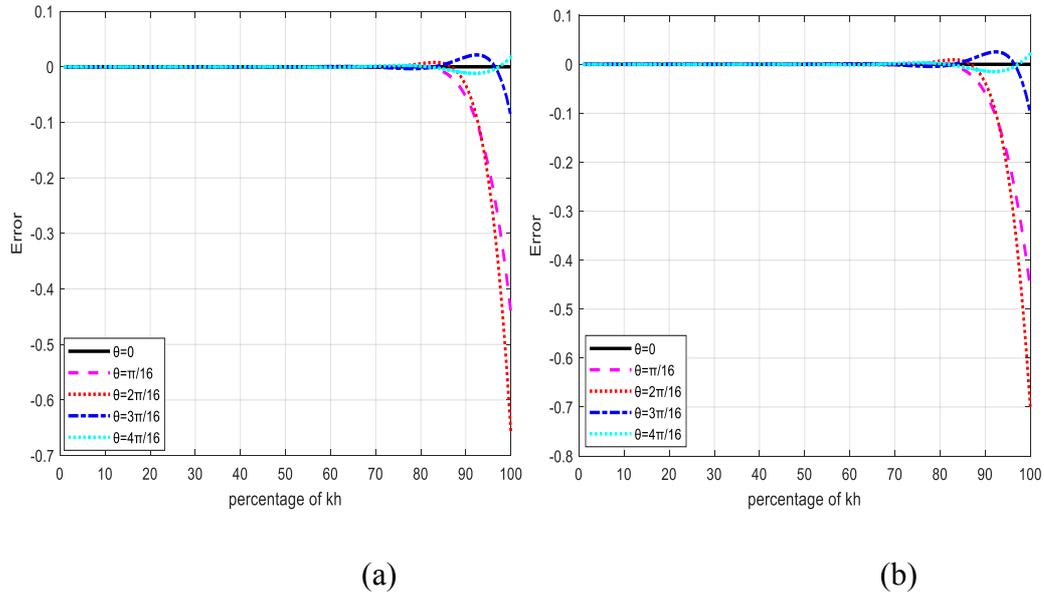

Figure 4. Dispersion error of the mixed derivative $\left(k_x k_z h^2\right)^2$ (a) The balanced SGFD scheme with *M*=7 and (b) the non-balanced SGFD scheme with *M*=7.

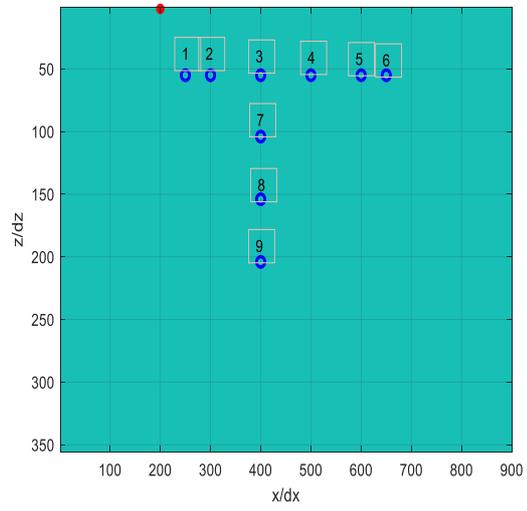

Figure 5. Source and receiver location

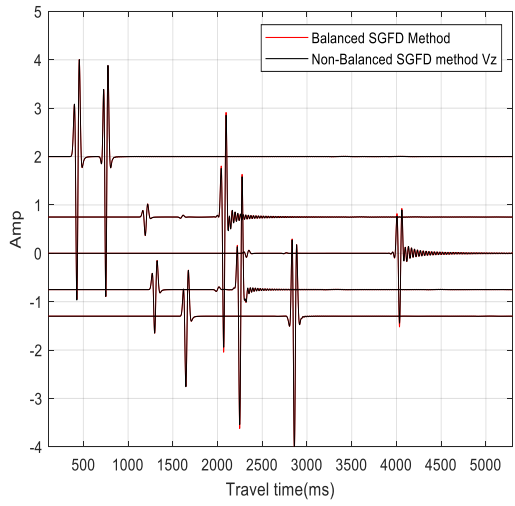 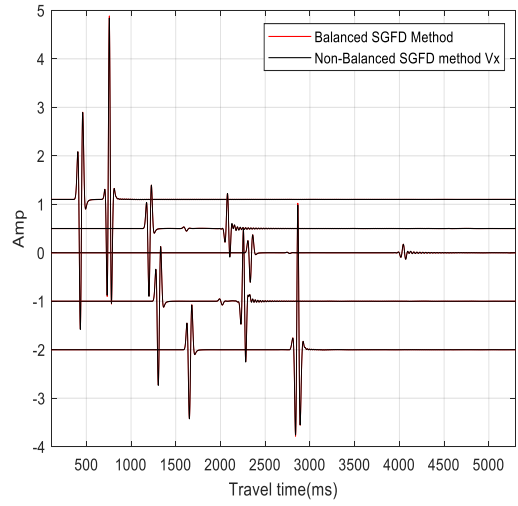

(a) (b)

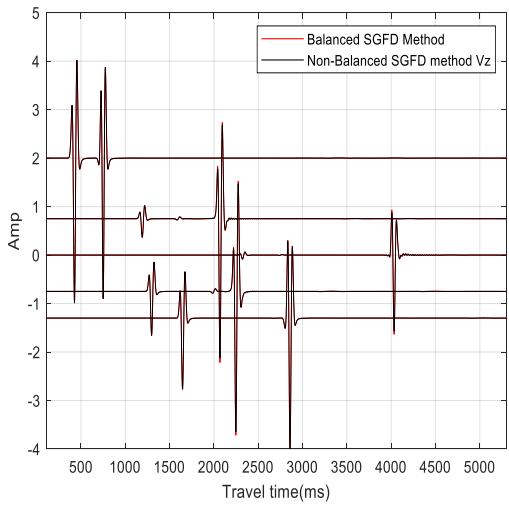 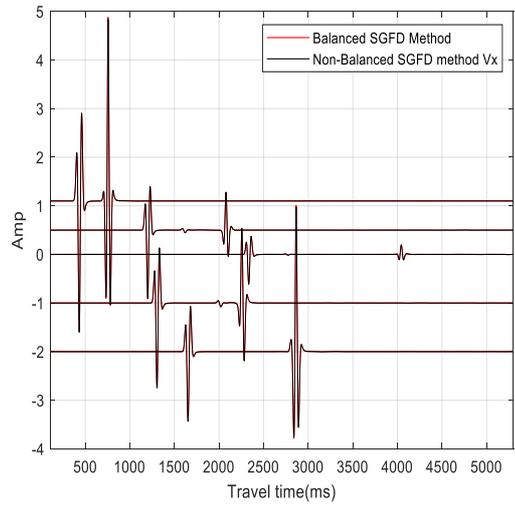

(c) (d)

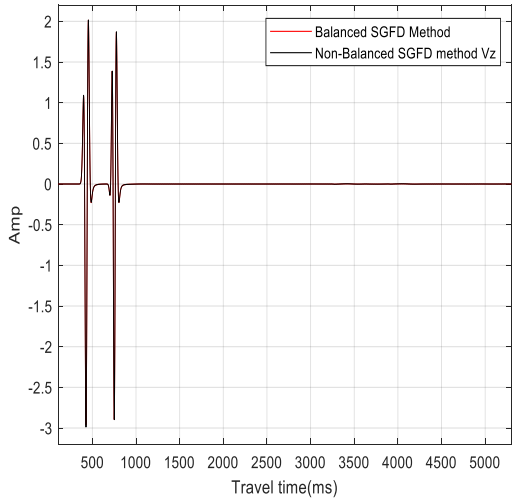
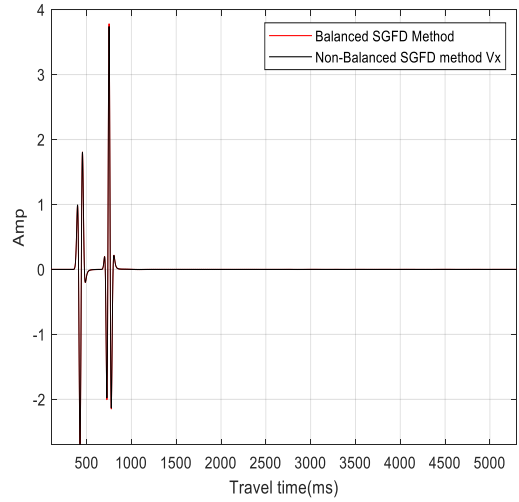

(e) (f)

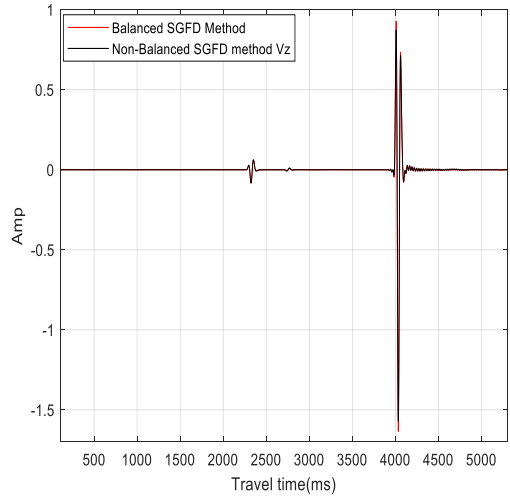
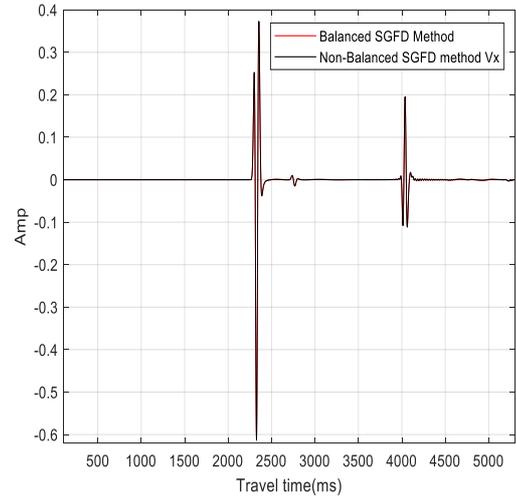

(g) (h)

Figure 6. Comparison of seismograms obtained at the 1,3,5,7,9 receivers from the top to the bottom (a) the $v_z$ component obtained with the balanced and the non-balanced SGFD scheme with $M$=4; (b) the $v_x$ component obtained with the balanced and the non-balanced SGFD scheme with $M$=4; (c) the $v_z$ component obtained with the balanced and the non-balanced SGFD scheme with *M*=7; (d) the $v_x$ component obtained with the balanced and the non-balanced SGFD scheme with *M*=7;(e) the local enlargement of the first seismogram in figure (c); (e) the local enlargement of the first seismogram in figure (d); (g) the local enlargement of the third seismogram in figure (c); (h) the local enlargement of the third seismogram in figure (d).

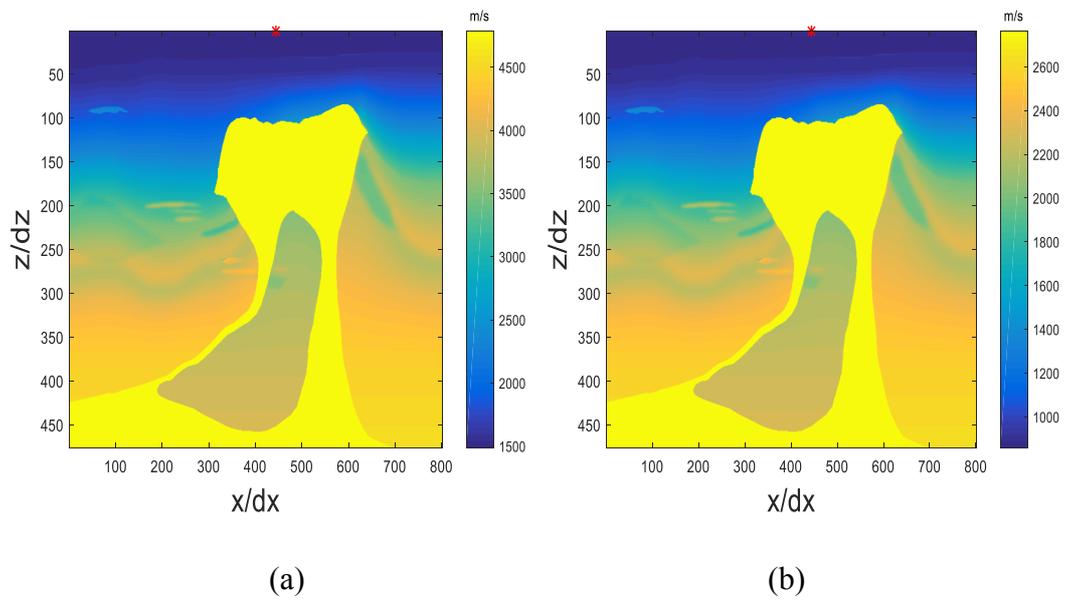

(a)                  (b)

Figure 7. Salt velocity model: (a) *P* wave velocity model; (b) *S* wave velocity model.

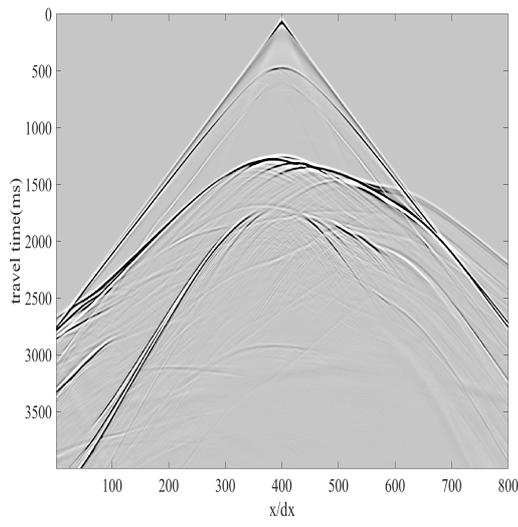 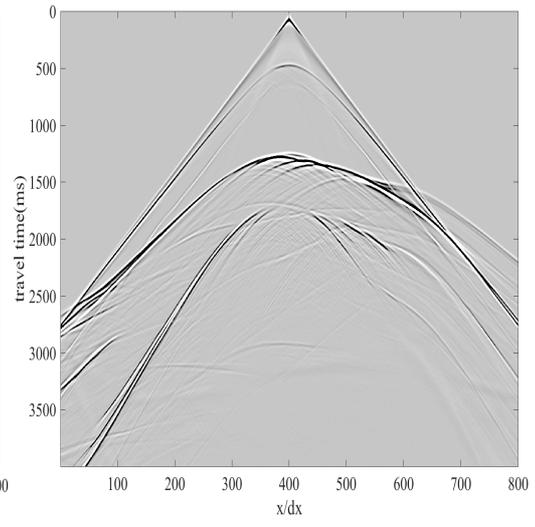

(a) (b)

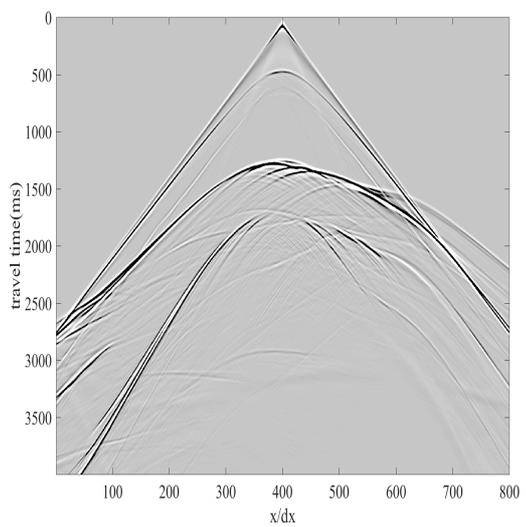 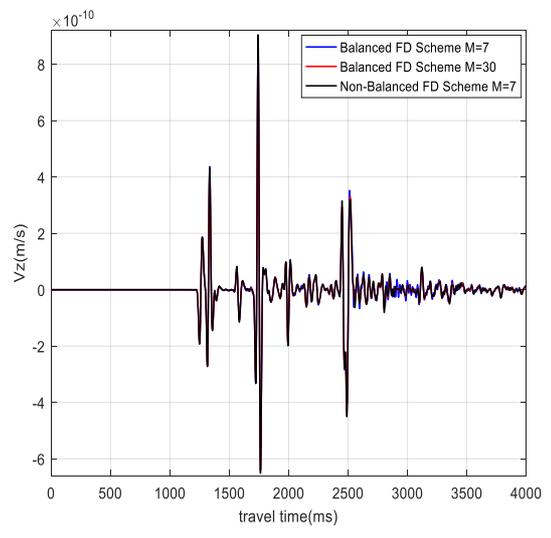

(c) (d)

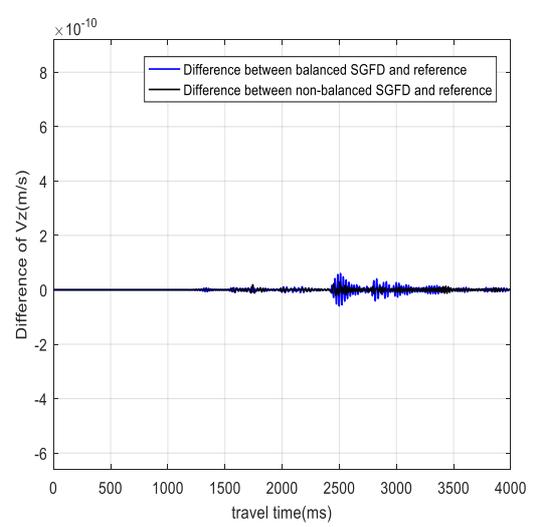

(e)

Figure 8. Seismic records of the $v_z$ component: (a) the balanced SGFD scheme (M=7); (b) the non-balanced SGFD scheme (M=7); (c) the balanced SGFD scheme (M=30); (d) slices of the seismograms extracted from (a)-(c) at x/dx=222; (e) difference of seismograms from (d).

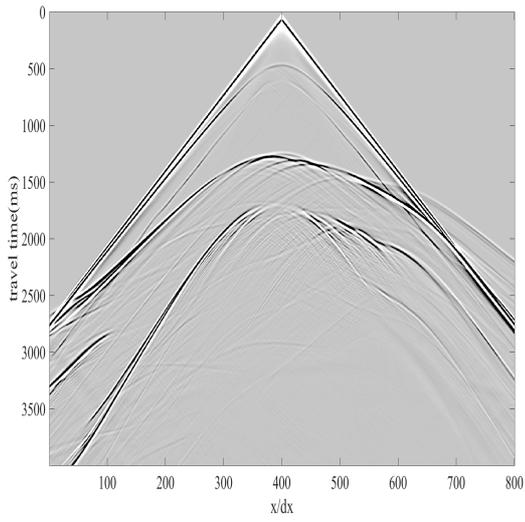 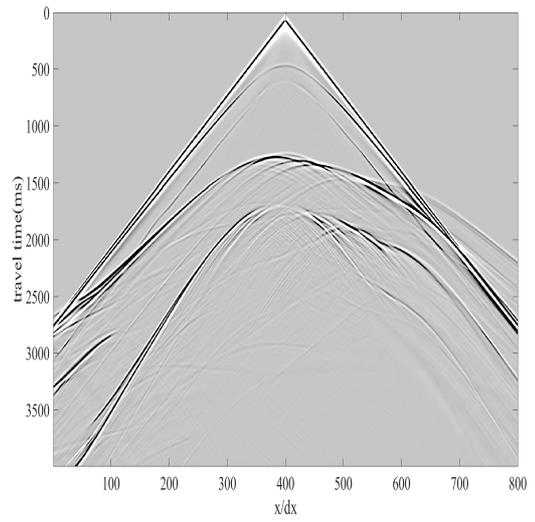

(a) (b)

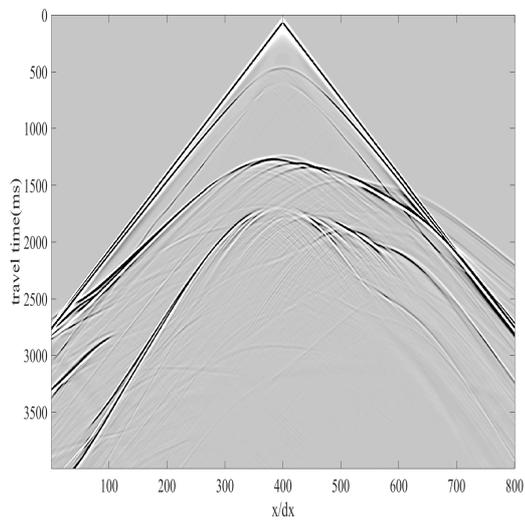 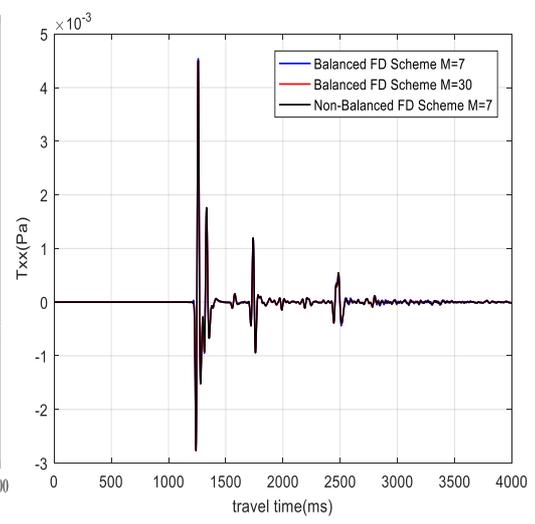

(c) (d)

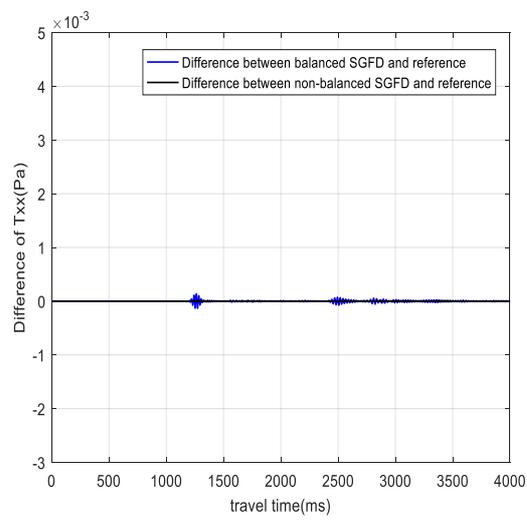

(e)

Figure 9. Seismic records of $\tau_{xx}$ component: (a) the balanced SGFD scheme; (b) the non-balanced SGFD scheme; (c) the balanced SGFD scheme (M=30); (d) slices of the seismograms extracted from (a)-(c) at x/dx=222; (e) difference of seismograms from (d).

Table 1. Optimized first-order SGFD coefficients for the non-balanced FD scheme.

|       | $C_1$   | $C_2$     | $C_3$     | $C_4$      | $C_5$      | $C_6$        | $C_7$       |
|-------|---------|-----------|-----------|------------|------------|--------------|-------------|
| $M = 3$ | 1.40887 | −0.16472  | 0.0172717 |            |            |              |             |
| $M = 5$ | 1.53147 | −0.252544 | 0.0607465 | −0.0135055 | 0.00199132 |              |             |
| $M = 7$ | 1.59906 | −0.310692 | 0.10345   | −0.0398274 | 0.0150857  | −0.00487876  | 0.00104241  |